\documentclass[a4paper]{article}

\usepackage{amsmath,amsthm,amsfonts}
\usepackage[dvips]{graphicx}

\newtheorem{thm}{Theorem}[section]
\newtheorem{defn}[thm]{Definition}
\newtheorem{lemma}[thm]{Lemma}
\newtheorem{prop}[thm]{Proposition}
\newtheorem{cor}[thm]{Corollary}
\newtheorem{que}[thm]{Question}

\theoremstyle{remark}   
\newtheorem{rem}[thm]{Remark}
\newtheorem{ex}[thm]{Example}

  \newcommand{\R}{\mathbb{R}}
 \newcommand{\T}{\mathbb{T}} \newcommand{\Z}{\mathbb{Z}}
\newcommand{\cA}{{\mathcal{A}}} \newcommand{\cC}{{\mathcal{C}}}
\newcommand{\cI}{{\mathcal{I}}} \newcommand{\cO}{{\mathcal{O}}}
\newcommand{\cL}{{\mathcal{L}}} \newcommand{\cM}{{\mathcal{M}}}
\newcommand{\cP}{{\mathcal{P}}} \newcommand{\cR}{{\mathcal{R}}}
\newcommand{\cS}{{\mathcal{S}}}

\newcommand{\sh}[1]{\text{sh}(#1)}  
\newcommand{\ssh}[1]{\text{sh}_0(#1)}   
\newcommand{\conv}[1]{\text{conv}(#1)}  

\begin{document}

\title{Boundary rigidity for Lagrangian submanifolds, non--removable intersections,
and Aubry--Mather theory }

\author{Gabriel P. Paternain \\
Department of Pure Mathematics and Mathematical Statistics\\
University of Cambridge\\
Cambridge CB3 0WB, England \and  Leonid Polterovich\footnote
 {Supported by the United States--Israel Binational
 Science Foundation grant 1999086}
\\School of Mathematical Sciences\\Tel Aviv
  University\\69978 Tel Aviv, Israel \and Karl Friedrich
  Siburg\footnote{Supported by a Heisenberg grant of the Deutsche
  Forschungsgemeinschaft
}\\Fakult\"at f\"ur Mathematik\\Ruhr--Universit\"at
  Bochum\\44780 Bochum, Germany}

\maketitle
\newpage

\tableofcontents
\newpage

\section{Introduction and results}
\label{intro}

In the present note we continue a theme which goes back to Arnold's seminal survey
''First steps in symplectic topology'' \cite{Arn}. A hypersurface in a cotangent
bundle is called optical if it bounds a fiberwise strictly convex domain. Likewise,
a Lagrangian submanifold is called optical if it lies in an optical hypersurface; a
particularly important class of examples is given by invariant tori in classical
mechanics. Arnold suggested to look at optical Lagrangian submanifolds from the
symplectic topology point of view. Arnold's suggestion inspired a number of results
in this direction (see, e.g., \cite{birkhoff,BP2}).

In this paper, we go a step further and establish a \emph{boundary rigidity}
phenomenon which, roughly speaking, can be formulated as follows. Certain Lagrangian
submanifolds lying in an optical hypersurface cannot be deformed into the domain
bounded by that hypersurface. Furthermore, even when boundary rigidity fails, we
often find another phenomenon called \emph{non--removable intersection}: the
intersection between the deformed Lagrangian submanifold and the hypersurface
contains some distinguished, dynamically relevant set. This observation links the
theory of symplectic intersections with modern aspects of dynamical systems.

Finally, we discuss Lagrangian submanifolds lying in the open domain
bounded by some optical hypersurface. Although these submanifolds
cannot be interpreted as invariant sets anymore, they still appear in
a number of interesting situations in geometry and dynamics.

\subsection{Preliminaries and basic notations}
\label{prelim}

Let $\theta: T^* X \to X $ be the cotangent bundle of a closed manifold $X$,
equipped with the canonical symplectic form $\omega= d\lambda$ where $\lambda$ is
the Liouville 1--form. We write $\cO$ for the zero section, and denote by
$\mathcal{L}$ the class of all Lagrangian submanifolds of $T^*X$ which are
Lagrangian isotopic to $\mathcal{O}$. Given $\Lambda \in \cL$, the natural
projection $\theta {\big |}_{\Lambda} : \Lambda \to X$ induces an isomorphism
between the cohomology groups $H^1(X,\R)$ and $H^1(\Lambda,\R)$. The preimage of
$[\lambda|_{\Lambda}]$ under this isomorphism is called the {\it Liouville class} of
$\Lambda$ and is denoted by $a_{\Lambda}\in H^1(X,\R)$. We say that a Lagrangian
submanifold $\Lambda \in \cL$ is {\it exact} if $a_{\Lambda} = 0$ and denote by
$\cL_0$ the class of all exact Lagrangian submanifolds in $\cL$.

A smooth, closed, fiberwise strictly convex hypersurface $\Sigma \subset T^*X$ is
called \emph{optical}. Fiberwise strict convexity means that $\Sigma$ intersects
each fiber $T_x^*X$ along a hypersurface whose second fundamental form is positive
definite. Denote by $\sigma$ the characteristic foliation of $\Sigma$, i.e., the
1--dimensional foliation tangent to the kernel of $\omega|_{T\Sigma}$. Note that
$\sigma$ is orientable and tangent to each Lagrangian submanifold contained in
$\Sigma$.

An orientable 1--dimensional foliation on a closed manifold is called {\it
conservative} if it admits a non--vanishing tangent vector field whose flow
preserves a measure which is absolutely continuous with respect to some (and hence
any) Riemannian measure on that manifold.

Let $\Lambda\in\cL$ be a Lagrangian submanifold lying in an optical hypersurface
$\Sigma$. Assume, in addition, that the restriction $\sigma|_{\Lambda}$ of the
characteristic foliation is conservative. In this case, one can show that $\Lambda$
is a \emph{section} of the cotangent bundle; this, a multidimensional version of the
Birkhoff second theorem, was established in \cite{birkhoff}. The assumption on the
conservativity of $\sigma|_{\Lambda}$ can be somewhat relaxed, but it is still
unknown whether it can be omitted completely. Interestingly enough, the same
assumption appears in a crucial way in the following, seemingly different context.

\subsection{Boundary rigidity}
\label{brp}

Suppose $\Sigma$ is a hypersurface bounding a domain $U_{\Sigma}$ and containing
some Lagrangian submanifold $\Lambda$. Can one push $\Lambda$ inside $U_{\Sigma}$ by
an exact Lagrangian deformation (i.e., a Lagrangian deformation preserving the
Liouville class)? We will present situations, both for the convex and nonconvex
case, where this is impossible. Sometimes, one cannot even move $\Lambda$ at all---a
phenomenon we call \emph{boundary rigidity}.

\subsubsection{The convex case}

Given an optical hypersurface $\Sigma$, we denote by $U_{\Sigma}$ the closed domain
bounded by $\Sigma$.

\begin{thm}
\label{brigid}
Let $\Lambda\in\cL$ be a Lagrangian submanifold lying in an optical hypersurface
$\Sigma$. Assume that the restriction $\sigma|_{\Lambda}$ of the characteristic
foliation is conservative. Let $K\in\cL$ be any Lagrangian submanifold\footnote{We
denote Lagrangian submanifolds by Greek letters, so this is a capital $\kappa$ and
not a capital $k$\ldots} lying in $U_{\Sigma}$ with the same Liouville class $a_K =
a_{\Lambda}$. Then $K=\Lambda$.
\end{thm}

In particular, $\Lambda$ cannot be deformed inside $U_{\Sigma}$ by an exact
Lagrangian isotopy, i.e., by a Lagrangian isotopy that preserves the Liouville
class.

Theorem~\ref{brigid} is proved in Section~\ref{selec} below. As the following
example shows, the assumption about the dynamics of the characteristic foliation
cannot be omitted.

\begin{ex}
\label{notrigid}
Consider $\Sigma = \{ H= 1 \} \subset T^* \T^2$ where
\begin{equation}
\label{hamex}
H(x_1,x_2,y_1,y_2) = (y_1-\sin x_1)^2 + (y_2-\cos x_1)^2.
\end{equation}
Then $\Sigma$ contains the zero section $\Lambda=\cO$. However, the restriction
$\sigma|_{\cO}$ of the characteristic foliation is a Reeb foliation with exactly two
limit cycles and, therefore, not conservative. We claim that $\cO$ is not boundary
rigid either. Indeed, the exact Lagrangian torus $K = \text{graph}(df)$ with
$$ f(x_1,x_2) = -\cos x_1 $$ does lie in $U_{\Sigma}$.
It is worth mentioning that $K$ intersects $\Sigma$ precisely
at the two limit cycles of the characteristic foliation. As
we will see in Section~\ref{nrcon}, this is no
coincidence.
\end{ex}

\subsubsection{The nonconvex case}
\label{nonconv1}

In this section, we let $\Sigma \subset T^*X$ be a smooth closed hypersurface which
need \emph{not} be strictly convex. Denote by $\sigma$ the characteristic foliation
of $\Sigma$.

We introduce the following class $\cL'$ of closed connected Lagrangian submanifolds
of $T^*X$. We say that $\Lambda \in \cL'$ if there exists a closed 1--form $\alpha$
on $X$ such that the restriction of $\lambda - \theta^*\alpha$ to $\Lambda$ is
exact, where $\theta: T^*X\to X$ is the natural projection. We call $a_{\Lambda}
:=[\alpha] \in H^1(X,\R)$ the Liouville class of $\Lambda$. Clearly, $\cL \subset
\cL'$, and the new notion of Liouville class agrees with the one given in
Section~\ref{prelim} for Lagrangian submanifolds in $\cL$.

Recall that a 1--dimensional foliation on a closed manifold is called {\it minimal}
if its leaves are everywhere dense. The following result is an immediate consequence
of Theorem~\ref{nonrem1} in Section~\ref{nrsect}.

\begin{thm}
\label{nonrem2}
Let $\Lambda\in \cL$ be a Lagrangian submanifold lying in $\Sigma$. Suppose that the
restriction $\sigma|_{\Lambda}$ of the characteristic foliation is minimal. Let $K
\in \cL'$ be any Lagrangian submanifold lying in $U_{\Sigma}$ with the same
Liouville class $a_K = a_{\Lambda}$. Then $K = \Lambda$.
\end{thm}

Comparing Theorem~\ref{nonrem2} to Theorem~\ref{brigid}, we see that the first is
applicable to a wider class of hypersurfaces (it does not require strict convexity);
on the other hand, the dynamical assumption on the characteristic foliation is more
restrictive.

\begin{ex}
\label{flat}
Consider the case when $\Sigma$ is the unit sphere bundle of the Euclidean metric on
the torus $\T^n$, i.e., $\Sigma = \{|p| = 1\}$. Theorem~\ref{brigid} yields boundary
rigidity for all flat Lagrangian tori $\{ p=v\}$ lying in $\Sigma$, while
Theorem~\ref{nonrem2} gives boundary rigidity only when the coordinates of $v$ are
rationally independent.
\end{ex}

\subsection{Non--removable intersections}
\label{nrint}

Let $\Sigma \subset T^* X$ be a hypersurface bounding a closed domain $U_{\Sigma}$.
Even if boundary rigidity, as explained in Section~\ref{brp}, may fail for a given
Lagrangian submanifold $\Lambda$, we will see that in many cases the intersection
$\Lambda \cap \Sigma$ cannot be empty. We call this phenomenon \emph{non--removable
intersection}. In fact, the intersection will always contain an invariant set of the
characteristic foliation of $\Sigma$.

\subsubsection{The nonconvex case}
\label{nrsect}

In this section, we let $\Sigma \subset T^*X$ be a smooth closed hypersurface which
need \emph{not} be strictly convex.

\begin{thm}
\label{nonrem1}
Let $\Lambda\in \cL$ be a Lagrangian submanifold lying in $\Sigma$. Let $K \in \cL'$
be any Lagrangian submanifold lying in $U_{\Sigma}$ with the same Liouville class
$a_K = a_{\Lambda}$. Then the intersection $K \cap \Lambda$ contains a compact
invariant set of the characteristic foliation $\sigma|_{\Lambda}$.
\end{thm}

Let us present the short and illustrative proof right here.

\begin{proof}
Let $v$ be a non--vanishing vector field on $\Sigma$ tangent to the characteristic
foliation $\sigma$. Arguing by contradiction, we suppose that the intersection $K
\cap \Lambda$ contains no compact invariant set of $\sigma$. Then, by a theorem of
Sullivan \cite{Su,LS2}, there exists a smooth function $h: \Sigma \to \R$ such that
$dh(x)\cdot v > 0$ for every $x \in K \cap \Lambda$. Extend $h$ to a function $H$
defined near $\Sigma$, and denote by $X_H$ the corresponding Hamiltonian vector
field. Note that $X_H$ is transverse to $\Sigma$ at each point in $K \cap \Lambda$.
Changing, if necessary, the sign of $v$ we can achieve that at these points $X_H$ is
pointing inside the domain $U_{\Sigma}$.

Let $\phi_t$ be the Hamiltonian flow of $X_H$, defined in a neighborhood of
$\Sigma$. Since $K\subset U_{\Sigma}$ and $\Lambda\subset \Sigma$, it follows that
$\phi_t(K) \cap \Lambda = \emptyset$, provided $t>0$ is sufficiently small. But the
Lagrangian submanifolds $\phi_t(K)$ and $\Lambda$ have equal Liouville classes, and
a Lagrangian intersection result due to Gromov\cite{gromov}---saying that two
Lagrangian submanifolds $\Lambda\in\cL$ and $K\in\cL'$ with equal Liouville classes
must intersect---gives the desired contradiction.
\end{proof}

\begin{rem}
Note that Theorem~\ref{nonrem1} immediately implies Theorem~\ref{nonrem2} about
boundary rigidity. Thus, boundary rigidity is a particular case of non--removable
intersections.
\end{rem}

In the following, it will be convenient to use the language of
symplectic shapes introduced by Sikorav \cite{sik,sik1} and Eliashberg
\cite{elshape}. The \textit{shape} of a subset $U \subset T^*X$ is
defined as
\begin{equation}
\label{defshape}
\sh{U} := \{a_{\Lambda} \in H^1(X,\R) \mid \Lambda\in \cL, \Lambda \subset U\}.
\end{equation}
Gromov's theorem \cite{gromov} implies Sikorav's elegant reformulation \cite{sik} of
Arnold's Lagrangian intersection conjecture proved in \cite{LS,H,gromov,Ch}: shapes
of disjoint subsets of $T^*X$ are disjoint. Therefore, every Lagrangian submanifold
$\Lambda\in \cL$ whose Liouville class lies on the boundary
$\partial\sh{U_{\Sigma}}$ must intersect $\Sigma$. Applying Sullivan's theorem
exactly as in the proof of Theorem~\ref{nonrem1}, one can refine this observation as
follows.

\begin{prop}
\label{nonrem11}
Suppose $\Lambda \in \cL$ is a Lagrangian submanifold such that $\Lambda\subset
U_{\Sigma}$ and $a_{\Lambda}\in \partial \sh{U_{\Sigma}}$. Then the intersection set
$\Lambda \cap \Sigma$ contains a compact invariant set of the characteristic
foliation of $\Sigma$.
\end{prop}

It would be interesting to understand, first, the dependence of the intersection set
$K \cap \Sigma$ on the Liouville class $a_K$ of $K$, and, secondly, the dynamical
meaning of the invariant set in the non--removable intersection. In fact, for the
convex case (i.e., when $\Sigma$ is optical) the latter can be done using
Aubry--Mather theory.

\subsubsection{The convex case}
\label{nrcon}

Suppose $\Sigma\subset T^*X$ is an optical hypersurface containing a Lagrangian
submanifold $\Lambda\in \cL_0$. As mentioned before, this setting allows a more
detailed description of non--removable intersections. Instead of going into
technicalities here, let us just illustrate this by an example.

\begin{ex}
Recall Example~\ref{notrigid} where $\Sigma=H^{-1}(1)$ with $H$ given by
\eqref{hamex}. The characteristic foliation has exactly two limit cycles, and we
constructed a particular Lagrangian submanifold $K$ such that $K\cap \Sigma$
consists of those two cycles. In fact, this is no coincidence, as the following
result shows.
\end{ex}

\begin{thm}
\label{twocycles}
Let $\Lambda\in\cL_0$ be any Lagrangian submanifold with $\Lambda\subset U_{\Sigma}$.
Then $\Lambda\cap\Sigma$ contains the two limit cycles of the characteristic
foliation of $\Sigma$.
\end{thm}

Precise results are stated and proved in Section~\ref{nonrem0}. The
methods we use are based on Aubry--Mather theory, which will be
explained in Section~\ref{abm}.

\subsection{Symplectic shapes of open fiberwise convex domains}
\label{shapes}

In this section, we focus on Lagrangian submanifolds lying in open, fiberwise convex
subsets of some cotangent bundle. Recall the definition of the symplectic shape for
a subset $U \subset T^*X$ given in \eqref{defshape}. In addition, we also define the
\emph{sectional shape}
$$ \ssh{U} \subset \sh{U} $$ to be the collection of all $a \in
H^1(X,\R)$ such that $U$ contains a Lagrangian section of $T^*X$ with Liouville
class $a$ (or, in other words, the graph of a closed 1--form representing the
cohomology class $a$). In contrast to the symplectic shape, $\ssh{U}$ is, in
general, \emph{not} preserved by Hamiltonian diffeomorphisms of $T^*X$ and, hence,
does not belong to the purely symplectic realm. However, it naturally arises and
plays a significant role in a number of interesting situations. In fact, this
antithesis was the starting point of our present research. It is resolved in a way
by the following theorem, which states that for open fiberwise convex sets $U
\subset T^*X$ both notions coincide.

\begin{thm}
\label{corf1}
Let $U\subset T^* X$ be an open fiberwise convex subset. Then every class $a\in
\sh{U}$ can be represented by a Lagrangian \emph{section} of the cotangent bundle.
In other words, $$ \ssh{U} = \sh{U}. $$
\end{thm}

By taking convex combinations, the following is a direct consequence of Theorem
\ref{corf1}.

\begin{cor}
\label{corf2}
The shape of an open fiberwise convex subset of $T^*X$ is an open convex subset of
$H^1(X,\R)$.
\end{cor}

Note that the shape of an open subset is always open (this follows immediately from
Weinstein's Lagrangian neighborhood theorem), so the main statement here is about
convexity. The proof of Theorem~\ref{corf1} is given in Section~\ref{2.3}.

\begin{ex}
Take a Riemannian metric $g$ on $X$ and consider the open unit ball bundle
$$U = \{(x,p) \in T^*X \mid |p|_g < 1\}.$$
There exists a remarkable norm on $H^1(X,\R)$, called {\it Gromov--Federer stable
norm}. Let us illustrate the corresponding dual norm $\Vert A \Vert$ for a homology
class $A \in H_1(X,\Z)$. Write $\ell (A)$ for the minimal length of a closed
geodesic representing $A$. Then
$$ \Vert A \Vert = \lim_{k \to \infty} \frac{\ell(kA)}{k} .$$

Gromov showed \cite{gromet} that the open unit ball of the stable norm coincides
with the sectional shape of $U$. In view of Theorem~\ref{corf1}, this is equal to
$\sh{U}$. Thus, for the Riemannian case, Theorem~\ref{corf1} leads to a geometric
description of the symplectic shape of a Riemannian unit ball bundle and, vice
versa, to a symplectic characterization of the unit stable norm ball.
\end{ex}

\begin{ex}
\label{exmather}
Let $H: T^*X \to \R$ be a fiberwise strictly convex Hamiltonian function. Assume
that $H$ has superlinear growth. Define the function $\alpha: H^1(X,\R) \to \R$ by
$$ \alpha (a) := \inf \{ h\in\R \mid a \in \ssh{\{H<h\}} . $$
This function is known as the convex conjugate of the \emph{Mather minimal action}
\cite{M1}; it was intensively studied in the past decade. Again, Theorem~\ref{corf1}
translates the variational definition of Mather's minimal action into symplectic
language.

As an illustration, consider the value $\min\alpha$. It is called \emph{Ma\~n\'e's
strict critical value}. It plays an important role when one studies the dependence
of the dynamics in the energy levels $\{H=h\}$ on the energy value $h$. This problem
is still far from being solved completely, even in a basic model of the magnetic
field on a closed manifold $X$ \cite{burns,pey}. It was proved in \cite{CIPP1} that
for $h > \min\alpha$, the dynamics in the energy level $\{H=h\}$ can be seen as a
time reparametrization of an appropriate Finsler flow on $X$. Ma\~n\'e's critical
value will appear in a crucial way in Sections~\ref{abm}--\ref{nonrem0} below.
\end{ex}

\section{Graph selectors of Lagrangian submanifolds and boundary rigidity}
\label{selec}

The main symplectic ingredient of our approach to proving Theorem~\ref{brigid} is
supplied by the following theorem which was outlined by Sikorav (in a talk held in
Chaperon's seminar) and proven by Chaperon (in the framework of generating
functions) and Oh (via Floer homology).

\begin{thm}[Sikorav, Chaperon \cite{chap}, Oh \cite{oh}]
\label{lip}
Let $\Lambda\subset T^* X$ be a Lagrangian submanifold in $\cL_0$. Then there exists
a Lipschitz continuous function $\Phi: X\to \R$, which is smooth on an open set
$X_0\subset X$ of full measure, such that
$$ (x,d\Phi(x)) \in \Lambda $$ for every $x\in X_0$. Moreover, if $d\Phi(x)
= 0$ for all $x\in X_0$ then $\Lambda$ coincides with the zero section $\cO$.
\end{thm}

We call the function $\Phi$ a \emph{graph selector} of the Lagrangian submanifold
$\Lambda$. In order to explain this terminology, consider $\Lambda$ as a
multi--valued section of the cotangent bundle. Then the differential $d\Phi(x)$
selects a single value of this section over the set $X_0$ in a smooth way.

In the following two sections, we will prove Theorem~\ref{lip} by using generating
functions quadratic at infinity, a powerful tool of symplectic topology in cotangent
bundles. Although this proof of Theorem~\ref{lip} is well known to experts, we were
unable to locate it in the literature.

In the final Section~\ref{rproof}, we apply Theorem~\ref{lip} in order to prove
Theorem~\ref{brigid}.

\subsection{Generating functions quadratic at infinity}
\label{gf}

Let $X$ be a closed manifold, and $E$ a finite dimensional real vector space. Denote
by $\cO_E$ the zero section of $T^*E$ and set
$$ V := T^*X \times \cO_E \subset T^*X \times T^*E = T^*(X \times E). $$

\begin{defn}
\label{defgf}
A smooth function $S: X \times E \to \R$ is called a \emph{generating function
quadratic at infinity (gfqi)} if $$ S(x,\xi) = Q_x(\xi)$$ outside a compact subset
of $X \times E$, where $Q_x$ is a smooth family of nondegenerate quadratic forms on
$E$, and graph$(dS)$ is transversal to $V$ in $T^*(X \times E)$.
\end{defn}

In particular, $W:= \text{graph}(dS) \cap V$ is a smooth closed submanifold of $V$
of the same dimension as $X$. Let $\chi: V \to T^*X$ be the natural projection. One
can show that the restriction of $\chi$ to $W$ is a Lagrangian immersion (see
\cite[Sect.~19]{AGV}). If $\chi|_W$ is an embedding then
$$ \Lambda:= \chi(W) $$ is a Lagrangian submanifold of $T^*X$. In this case
we say that $S$ is a gfqi of $\Lambda$; this means that
\begin{equation}
\label{eqgf}
\Lambda = \{(x,d_x S(x,\xi)) \mid x \in X, \xi \in E, d_{\xi}S(x,\xi) = 0 \}.
\end{equation}

Regarding the existence of a gfqi of a given Lagrangian submanifold $\Lambda$, it is
known that every $\Lambda\in\cL_0$ admits a gfqi \cite{sik0, Ch}. We do not know any
existence result for Lagrangian submanifolds $\Lambda\notin \cL_0$.

\subsection{The graph selector---proof of Theorem~\ref{lip}}

Let $S:X \times E \to \R$ be a gfqi of a Lagrangian submanifold $\Lambda \in \cL$.
The graph selector is defined by a suitable minimax procedure which we are going to
describe now.

Fix a scalar product on $E$. Let $B_x:E \to E$ be a self--adjoint
operator so that $Q_x(\xi) = (B_x \xi,\xi)$. Denote by $E^-_x$ the
subspace of $E$ generated by all eigenvectors of $B_x$ with negative
eigenvalues. Set
$$ E^a_x =\{\xi \in E \mid S_x(\xi) \leq a \} $$ where $a \in \R$ and
$S_x(\cdot) := S(x,\cdot)$. Take $N > 0$ so that $S(x,\xi) = Q_x(\xi)$
whenever $|Q_x(\xi)| \geq N$. The quadratic forms $Q_x,\; x \in X$
have the same index which we denote by $m$. The homology group
$H_m(E^N_x,E^{-N}_x;\Z_2)$ is isomorphic to $\Z_2$, and the generator,
say $A_x$, is represented by the $m$-dimensional disc in $E^-_{x}$
whose boundary lies in $\{Q_x(\xi) = -N\}$. For $a \in [-N,N]$,
consider the natural morphism
$$ I_{a,x} : H_m(E^a_x,E^{-N}_x;\Z_2) \to H_m(E^N_x,E^{-N}_x;\Z_2). $$

\begin{defn}
The function $\Phi: X\to\R$ defined by
$$ \Phi(x) := \inf\{a \mid A_x \in \text{Image}(I_{a,x}) \} $$ is
called the \emph{graph selector} of $\Lambda$ associated to the gfqi $S$.
\end{defn}

We claim that $\Phi$ has the properties stated in Theorem \ref{lip}. Clearly,
$\Phi(x)$ is a critical value of $S_x$. Consider the subset $X_0 \subset X$
consisting of all those $x$ for which $S_x$ is a Morse function whose critical
points have pairwise distinct critical values.

In a neighborhood $U$ of any point of $X_0$ there exists a smooth function $\varphi:
U \to E$ such that $\varphi(x)$ is a critical point of $S_x$ and $\Phi(x) =
S(x,\varphi(x))$. Differentiating with respect to $x$ and taking into account that
$d_{\xi}S(x,\varphi(x)) = 0$ we get that $d\Phi (x) = d_x S (x,\varphi(x))$. Thus,
in view of (\ref{eqgf}), $(x, d\Phi(x) ) \in \Lambda$ for all $x \in X_0$, so $\Phi$
is indeed a selector.

\begin{prop}
\label{x0}
$X_0$  is an open subset of $X$ of full measure.
\end{prop}

\begin{proof}
Let $\theta : T^*X \to X$ be the natural projection. A simple local analysis shows
that $S_x$ is Morse if and only if $x$ is a regular value of $\theta|_{\Lambda}$
(see \cite{AGV}, section 21.2). Denote the set of such $x\in X$ by $X_1$. It is an
open subset of $X$, and by Sard's Theorem it has full measure.

Let $U \subset X_1$ be a sufficiently small open subset. The critical
points of $S_x$ depend smoothly on $x \in U$. Denote them by
$\varphi_1(x),\ldots, \varphi_d(x)$, and put
$$a_{ij}(x) = S(x,\varphi_i(x)) - S(x,\varphi_j(x)) $$ for $i \neq
j$. Note that
$$da_{ij} (x) = d_x S(x,\varphi_i(x)) - d_x S(x,\varphi_j(x)) \neq 0$$
since the map
$$\chi|_W : W \to T^*X,\; (x,\xi) \mapsto (x, d_x S (x,\xi))$$ is an
embedding. Therefore the sets $\Gamma_{ij} = \{x \in U \mid a_{ij}(x)
= 0\}$ are smooth hypersurfaces. But, by definition of $X_0$, we have
$$X_0 \cap U = U \setminus \bigcup_{i\neq j} \Gamma_{ij}, $$ so $X_0
\cap U$ is an open subset of full measure in $X\cap U$.
\end{proof}

\begin{prop}
\label{un}
If $d\Phi(x) = 0$ for all $x \in X_0$ the submanifold $\Lambda$ coincides with the
zero section of $T^*X$.
\end{prop}

\begin{proof}
Identify $X$ with the zero section of $T^*X$.  Since $X_0$ has full measure, its
closure equals $X$. Hence $\Lambda$ contains $X$ since $d\Phi(x) = 0$ for $x \in
X_0$, and thus $\Lambda = X$.
\end{proof}

\begin{prop}
\label{lip1}
$\Phi$ is a Lipschitz function on $X$.
\end{prop}

\begin{proof}
Since $X$ is compact it suffices to prove this locally. Let $U \subset X$ be a
sufficiently small open subset. There exists a smooth family of linear automorphisms
$F_x: E \to E, x \in U$, and a quadratic form $Q$ on $E$, so that $Q_x \circ F_x =
Q$ for all $x \in U$. It is easy to see that the function $S'(x,\xi) := S(x,F_x
\xi)$ is again a gfqi of $\Lambda$ over $U$, whose graph selector coincides with
$\Phi|_U$. In what follows we work with $S'$ instead of $S$, because the functions
$S'_x, x \in U$, equal the same quadratic form $Q$ outside a compact subset of $E$.
Therefore there exists a positive constant $C$ such that for all $x,y \in U$ and
$\xi \in E$ we have
\begin{equation}
\label{lip2}
|S'(x,\xi) -S'(y,\xi)| \leq C|x-y|.
\end{equation}
Fix $\epsilon > 0$ and $x \in U$, and set
$$ a(y) := \Phi(x) + \epsilon + C|x-y|, $$ for all $y\in U$. It
follows from inequality \eqref{lip2} that $E^{a(x)}_x \subset
E^{a(y)}_y$ for all $y \in U$. By definition, the pair $(E^{a(x)}_x,
E^{-N})$ contains a relative cycle representing the class
$A_x$. Therefore, the same holds for the pair $(E^{a(y)}_y,
E^{-N})$. We get that $\Phi(y) \leq a(y)$, which yields
$$ \Phi (y) - \Phi(x) \leq C|x-y| + \epsilon. $$ Since the last
inequality is valid for each $\epsilon > 0$ we have
$$ \Phi (y) - \Phi (x) \leq C|x-y|. $$ Finally, interchanging $x$ and
$y$ we get that $\Phi$ is Lipschitz continuous.
\end{proof}

Thus, the function $\Phi$ satisfies all requirements of a graph selector. This
finishes the proof of Theorem \ref{lip}.
\hfill\qed

\subsection{Proof of Theorem \ref{brigid}}
\label{rproof}

Let $\Lambda\in\cL$ be a Lagrangian submanifold lying in some optical hypersurface
$\Sigma$, and assume that the restriction $\sigma|_{\Lambda}$ of the characteristic
foliation is conservative. Let $K\in\cL$ be any Lagrangian submanifold lying in
$U_{\Sigma}$ with the same Liouville class. We want to prove that $K=\Lambda$.

By the multidimensional Birkhoff theorem \cite{birkhoff}, $\Lambda$ is a Lagrangian
section, i.e., $\Lambda = \text{graph}(\alpha)$ for some closed 1--form $\alpha$. By
applying the symplectic shift $(x,p)\mapsto (x,p-\alpha(x))$ we may assume that
$\Lambda = \cO$ is the zero section. Note that the transformed hypersurface remains
optical.

Suppose now there is another Lagrangian submanifold $K \subset U_{\Sigma}$, obtained
from $\Lambda$ by an exact Lagrangian deformation. Let $ \Phi : X \to \R $ be a
graph selector of $K$ so that $(x,d\Phi(x)) \in K$ for all $x \in X_0$, where
$X_0\subset X$ is a set of full measure as in Theorem~\ref{lip}.

Pick a smooth Hamiltonian function $H: T^* X \to \R$ which is fiberwise strictly
convex such that $\Sigma$ is a regular level set of $H$. Since $\Lambda=\cO$ the
vector $\frac{\partial H}{\partial p}(x,0)$ gives the outer normal direction to the
hypersurface $\Sigma \cap T^*_x X\subset T^*_x X$. Because $K\subset U_{\Sigma}$, we
have
\begin{equation}
\label{ineq}
d\Phi(x)\cdot \frac{\partial H}{\partial p}(x,0) < 0
\end{equation}
in local canonical coordinates $(x,p)$ for all $x\in X_0$
with $d\Phi(x) \neq 0$.

Let $v$ be a non--singular vector field on $\Lambda$ which is tangent to the
characteristic foliation, and whose flow $\psi_s$ preserves a measure $\mu$ which is
absolutely continuous with respect to some Riemannian measure. Then the Hamiltonian
differential equations for $H$ show that $v$ is collinear to the vector field
$\frac{\partial H}{\partial p}(x,0) $ on $\Lambda$. In view of \eqref{ineq}, we may
assume that
\begin{equation}
\label{ineq1}
d\Phi(x)\cdot v(x) < 0
\end{equation}
for all $x\in X_0$ with $d\Phi(x) \neq 0$.

On the other hand, we claim that
\begin{equation}
\label{integral}
\int_{X_0} d\Phi(x)\cdot v(x) d\mu(x) = 0.
\end{equation}
Note that the theorem is an immediate consequence of \eqref{integral}. Indeed,
combining \eqref{integral} with \eqref{ineq1} we see that $d\Phi$ must vanish on
$X_0$, and hence $$ K = \cO=\Lambda
$$ in view of Theorem \ref{lip}.

It remains to prove formula \eqref{integral}. Since the
function $\Phi$ is Lipschitz continuous, the function
$s\mapsto \Phi(\psi_s x ) -\Phi (x)$ on $[0,1]$ is also
Lipschitz continuous for every $x \in X$. By Rademacher's
theorem, it is differentiable almost everywhere with
$$ \Phi(\psi_1 x ) -\Phi (x) = \int_0^1 \frac{d}{ds} \Phi (\psi_s x)
ds. $$ Since the flow $\psi_s$ preserves the measure $\mu$ we have
$$ 0 = \int_X [\Phi(\psi_1 x ) -\Phi (x)] d\mu(x) = \int_X \int_0^1
\frac{d}{ds} \Phi (\psi_s x) ds\, d\mu(x) .$$ Since $X_0$ has full
measure with respect to $\mu$
and since $\psi_s$
preserves $\mu$, we have
\begin{align*}
0 &= \int_0^1 \int_X \frac{d}{ds} \Phi(\psi_s x) d\mu(x)\,ds\\ &= \int_0^1
\int_{\psi_s^{-1}(X_0)} d\Phi(\psi_s x)\cdot v(\psi_s x) d\mu(x)\,ds\\ &= \int_0^1
\int_{X_0} d\Phi(x)\cdot v(x) d\mu(x)\,ds\\ &= \int_{X_0} d\Phi(x)\cdot v(x)
d\mu(x).
\end{align*}
This proves \eqref{integral} and finishes the proof of the theorem.
\hfill\qed

\section{Brief summary of Aubry--Mather theory}
\label{abm}

In this section, we give a brief overview about what is known as
Aubry--Mather theory. We refer the reader to the books \cite{CI, F}
for various preliminaries related to the material presented here.

\subsection{Ma\~n\'e's critical value}
\label{maneval}

Let $X$ be a closed connected smooth manifold and let $L:TX\to \R$ be a smooth,
fiberwise convex, superlinear Lagrangian\footnote{We always distinguish between the
term ``Lagrangian'' (i.e., Lagrangian function) and ``Lagrangian submanifold''.}.
This means that $L$ restricted to each $T_{x}X$ has positive definite Hessian and
for some Riemannian metric we have
$$ \lim_{|v|\to \infty}\frac{L(x,v)}{|v|}=\infty $$ uniformly on
$x\in X.$ Let $H:T^{*}X\to \R$ be the Hamiltonian associated to $L$ and
$$ \ell : TX \to T^*X $$ be the Legendre transform $\ell: (x,v)\mapsto
\frac{\partial L(x,v)}{\partial v}$. Since $X$ is compact, the extremals of $L$ give
rise to a complete flow $\phi_t$ on $TX$, called the Euler--Lagrange flow of the
Lagrangian. Using the Legendre transform we can push forward $\phi_t$ to obtain
another flow $\phi^*_t$ on $T^*X$ which is the Hamiltonian flow of $H$ with respect
to the canonical symplectic structure of $T^*X$. The energy of $L$ is the function
$E:TX\to \R$ given by
$$ E(x,v) := \frac{\partial{L}}{\partial{v}}(x,v)\cdot v-L(x,v) = H(\ell(x,v)). $$
The energy $E$ is a first integral of the Euler--Lagrange flow $\phi_{t}$.

Recall that the $L$--action of an absolutely continuous\footnote{A curve $\gamma:
[a,b] \to X$ is called absolutely continuous if for every $\epsilon > 0$ there
exists $\delta > 0$ so that for each finite collection of pairwise disjoint open
intervals $(s_i,t_i)$ in $[a,b]$ of total length $< \delta$ one has $\sum_{i=1}^N
\text{dist}(\gamma(t_i),\gamma(s_i)) < \epsilon$. Here $\text{dist}$ is any
Riemannian distance on $X$.} curve $\gamma:[a,b]\rightarrow X$ is defined by
$$ A_{L}(\gamma) := \int_{a}^{b}L(\gamma(t),\dot{\gamma}(t))\,dt. $$
Given two points $x_{1},x_{2}\in X$ and some $T>0$ denote by $\cC_{T}(x_{1},x_{2})$
the set of absolutely continuous curves $\gamma:[0,T]\rightarrow M$ with
$\gamma(0)=x_{1}$ and $\gamma(T)=x_{2}$. For each $k\in \R$, we define
$$ \Phi_{k}(x_{1},x_{2};T):=\inf\{A_{L+k}(\gamma) \mid
\gamma\in\cC_{T}(x_{1},x_{2})\}. $$ The {\it action potential} $\Phi_k: X\times X\to
\R\cup \{-\infty\}$ of $L$ is defined by
$$ \Phi_{k}(x_{1},x_{2}) := \inf_{T>0}\Phi_{k}(x_{1},x_{2};T). $$

\begin{defn}[Ma\~n\'e]
\label{defmane}
The \emph{critical value} of $L$ is the real number
$$ c=c(L) := \inf \{ k\in\R \mid \Phi_{k}(x,x)>-\infty \text{ for
some } x\in X \}.
$$
\end{defn}

Note that actually $\Phi_{k}(x,x)>-\infty$ for \emph{all} $x\in X$. Since $L$ is
convex and superlinear, and $X$ is compact, such a number exists. It singles out the
energy level in which relevant globally action--minimizing orbits and/or measures
live \cite{DC,Ma2,CDI,F}. Their study has a long history that goes back M. Morse and
G.A. Hedlund; recently, there has been a great deal of activity on this subject,
cf.\ \cite{BP0,B1,B2,burns,DC,F,Ma2,Ma3,M1,M2,pey,sib}.

The critical value can be characterized in a variety of ways
\cite{Ma2,CDI,CIPP1,CIPP2}. Each of these characterizations gives a new insight into
geometry and dynamics. Let us explain first the relation of the critical value with
Mather's theory of minimizing measures \cite{M1}.

Let $\cP(L)$ be the set of Borel probability measures on $TX$ that have compact
support and are invariant under the Euler--Lagrange flow $\phi_{t}$. Ma\~n\'e
\cite{Ma2,CDI} showed that the critical value can be described in terms of measures
as
\begin{equation}
\label{ergodic}
c(L)=-\min\left\{\int L\,d\mu \mid \mu\in{\mathcal P}(L)\right\}.
\end{equation}
We will say that $\mu\in{\mathcal P}(L)$ is a \emph{minimizing measure} if $\mu$
realizes the minimum in (\ref{ergodic}). The \emph{Mather set} in $TX$ is defined as
$$ \tilde\cM := \overline{\bigcup_{\mu}\text{supp}(\mu)}, $$ where
supp($\mu$) is the support of the measure $\mu$, the bar denotes the closure of a
set, and the union is taken over all minimizing measures. Mather's Lipschitz graph
theorem \cite{M1} asserts that $\tilde{{\cM}}$ is a Lipschitz graph with respect to
the canonical projection $\tau:TX\to X$. We call $\cM := \tau(\tilde\cM)\subset X$
the \emph{projected Mather set}. It is known that $\tilde{{\cM}}$ is contained in
the energy level $E^{-1}(c)$ \cite{DC}. We define the \emph{Mather set}
$\tilde\cM^*$ in $T^*X$ as the image of $\tilde\cM$ under the Legendre transform.

It turns out the critical value $c(L)$ can be recovered purely from the Hamiltonian
as the following result obtained in \cite{CIPP1} (and also independently by Fathi)
shows. Namely, we have
\begin{equation}
\label{minimax}
c = c(H) = \inf_{u\in C^{\infty}(X,\R)}\max_{x\in X} H(x,du(x)).
\end{equation}

In fact, Theorem \ref{corf1} gives a new, more geometric way of
looking at this quantity (cf.\ Example \ref{exmather}). It implies that
\begin{equation}
\label{newway}
c = \inf_{\Lambda\in \cL_0} \max_{(x,p)\in \Lambda} H(x,p)
\end{equation}
where, as usual, $\cL_0$ denotes the class of exact Lagrangian submanifolds in
$\cL$.

\subsection{Weak KAM solutions and Peierls barrier}
\label{wkamsol}

Given a continuous function $u:X\to\R$, we write
$$ u \prec L+c $$ whenever $u(x)-u(y)\le\Phi_c(y,x)$ for all $x,y\in
X$. Here, $\Phi_c$ is the action potential for the critical value $c$.

\begin{rem}
\label{smaller}
The condition $u\prec L+c$ is actually equivalent to $u$ being
Lipschitz and $H(x,du(x))\leq c$ for almost every $x\in X$
\cite[Thm.~4.2.10\& Lemma~4.2.11]{F}. Recall that by Rademacher's
theorem, Lipschitz functions are differentiable almost everywhere.
\end{rem}

We say that a continuous function $u_{+}:X\to\R$ is a {\it positive weak KAM
solution} if $u_{+}$ satisfies the following two conditions:
\begin{enumerate}
\item $u_{+}\prec L+c$; \item for all $x\in X$, there exists a absolutely continuous
curve $\gamma_{+}^{x}:[0,\infty)\to X$ such that $\gamma_{+}^{x}(0) = x$ and
$$ u_{+}(\gamma_{+}^{x}(t))-u_{+}(x)=
\int_{0}^{t}(L+c)(\gamma_{+}^{x}(s),\dot{\gamma}_{+}^{x}(s))\,ds $$ for all $t\geq
0$.
\end{enumerate}
Similarly we say that a continuous function
 $u_{-}:X\to\R$ is a {\it negative
weak KAM solution} if $u_{-}$ satisfies the following two conditions:
\begin{enumerate}
\item $u_{-}\prec L+c$; \item for all $x\in X$, there exists an absolutely
continuous curve $\gamma_{-}^{x}:(-\infty,0]\to X$ such that $\gamma_{-}^{x}(0) = x$
and
$$ u_{-}(x)-u_{-}(\gamma_{-}^{x}(-t))=
\int_{-t}^{0}(L+c)(\gamma_{-}^{x}(s),\dot{\gamma}_{-}^{x}(s))\,ds $$ for all $t\geq
0$.
\end{enumerate}

Fathi's weak KAM theorem asserts that positive and negative weak KAM solutions
always exist \cite{F}. At any point $x$ of differentiability of a weak KAM solution
$u$, Conditions 1 and 2 imply that
$$ H(x,du(x))=c. $$
In fact, the points $x$ of differentiablity of $u_{+}$ (resp. $u_{-}$) are precisely
those for which the curve $\gamma_{+}^{x}$ (resp. $\gamma_{-}^{x}$) is unique.

We denote by $\cS_{\pm}$ the set of all positive (respectively, negative) weak KAM
solutions. A pair of functions $(u_{-},u_{+})$ is said to be \emph{conjugate} if
$u_{\pm}\in {\cS}_{\pm}$ and $u_{-}=u_{+}$ on the projected Mather set $\cM$. We
will need the following result.

\begin{thm}[{\cite[Thm 5.1.2]{F}}]
\label{conjugate}
If $u:X\to \R$ is a function such that $u\prec L+c$, then there exists a unique pair
of conjugate functions $(u_{-},u_{+})$ such that $u_{+}\leq u\leq u_{-}$.
\end{thm}

The \emph{Peierls barrier} \cite{M2} is the function $h:X\times
X\to\R$ defined by
$$ h(x,y) := \liminf_{T\to \infty} \Phi_{c}(x,y;T). $$ The function
$h$ is Lipschitz (cf. \cite[Corollary 5.3.3]{F}) and, obviously,
satisfies $h(x,y)\geq \Phi_{c}(x,y)$. The Peierls barrier can be
recovered from the weak KAM solutions or the action potential
(cf. \cite[Prop.~13]{CIPP2}). Corollary 5.3.7 in \cite{F} gives
\begin{equation}
\label{peierls}
h(x,y) = \max_{(u_{-},u_{+})} (u_{-}(y)-u_{+}(x))
\end{equation}
where the maximum is taken over all pairs $(u_-,u_+)$ of conjugate functions.

\subsection{The Aubry set}
\label{aset}

By definition, two conjugate functions $u_{\pm}$ coincide on the projected Mather
set $\cM$. It turns out, however, that in general there is a bigger set, called
Aubry set, with this property. Namely, setting
$$ \cI_{(u_{-},u_{+})} := \{ x\in X \mid u_-(x)=u_+(x) \}, $$
we can define
$$ \cA := \bigcap_{(u_{-},u_{+})} \cI_{(u_{-},u_{+})} $$
where the intersection is taken over all pairs of conjugate functions. This set is
the \emph{projected Aubry set}. Clearly, it contains the projected Mather set $\cM$.

In order to define the Aubry set in $T^*X$, we note that the functions $u_{-}$ and
$u_{+}$ are differentiable at every point $x\in \cI_{(u_{-},u_{+})}$ with the same
derivative. Moreover, the map $\cI_{(u_{-},u_{+})}\ni x\mapsto du_-(x)= du_+(x) \in
T^*X$ is Lipschitz continuous. This was proved by Fathi \cite[Thm.~5.2.2]{F}. That
map defines a set
$$ \tilde\cI_{(u_-,u_+)} \subset T^*X $$ that projects injectively onto
$\cI_{(u_{-},u_{+})}$ and contains the Mather set. The \emph{Aubry set} in $T^*X$ is
defined as
$$ \tilde\cA^* := \bigcap_{(u_{-},u_{+})} \tilde\cI_{(u_-,u_+)}, $$ where, again, the
intersection is taken over all pairs $(u_{-},u_{+})$ of conjugate functions. It
turns out that $\cA =\theta(\tilde\cA^*)$. As usual, we denote the preimage of
$\tilde\cA^*$ under the Legendre transform by $\tilde\cA$ and call it the Aubry set
in $TX$. The sets $\tilde\cM$ and $\tilde\cA$ are compact and invariant under the
Euler--Lagrange flow $\phi_{t}$.

It turns out that the Aubry set consists of a distinguished kind of orbits. To make
this precise, we say that an absolutely continuous curve $\gamma:[a,b]\to X$ is {\it
semistatic} if
$$ A_{L+c}\left(\gamma\vert_{[s,t]}\right) =
\Phi_c(\gamma(s),\gamma(t)) $$ for all $a\le s\le t\le b$.  Semistatic curves are
solutions of the Euler--Lagrange equation because of their minimizing properties.
Also it is not hard to check that semistatic curves have energy precisely $c$
\cite{Ma2, CDI}. We say that an absolutely continuous curve $\gamma:[a,b]\to X$ is
{\it static} if it is semistatic and
$$ \Phi_c(\gamma(s),\gamma(t))+\Phi_c(\gamma(t),\gamma(s))=0 $$ for
all $a\le s\le t\le b$. The notions of static and semistatic curves are closely
related to Mather's notions of $c$--minimal trajectories and regular $c$--minimal
trajectories \cite{M2}.

\begin{prop}
\label{aubdes}
The Aubry set $\tilde\cA$ consists precisely of those orbits whose projections to
$X$ are static curves.
\end{prop}

This is well known to experts; nevertheless, we include its proof for the sake of
completeness.

\begin{proof}
Take $(x,v)\in \tilde\cA$. We show that $\gamma(t)=\tau(\phi_{t}(x,v))$ is a static
curve. By the definition of $\tilde\cA$ and Theorem 5.2.2 in \cite{F}, we have for
any pair $(u_{-},u_{+})$ of conjugate functions that
$$ u_{+}(\gamma(t))-u_{-}(\gamma(s))=
u_{+}(\gamma(t))-u_{+}(\gamma(s))=A_{L+c}\left(\gamma|_{[s,t]}\right)
$$ for all $s\leq t$. Using \eqref{peierls} we can choose a pair
$(u_{-},u_{+})$ of conjugate functions for which the Peierls barrier
$h$ satisfies
$$ h(\gamma(t),\gamma(s))=u_{-}(\gamma(s))-u_{+}(\gamma(t)). $$
Therefore, we can estimate
\begin{equation}
\label{aaa}
A_{L+c}(\gamma|_{[s,t]})+\Phi_{c}(\gamma(t),\gamma(s)) \leq
u_{+}(\gamma(t))-u_{-}(\gamma(s))+h(\gamma(t),\gamma(s))=0.
\end{equation}
It is easy to show that $\Phi_{c}$ satisfies the triangle inequality
\[ \Phi_c(x,y)\leq \Phi_c(x,z)+\Phi_c(z,y) \]
for all $x,y,z\in X$, as well as $\Phi_{c}(x,x)=0$ for all $x\in X$. Hence we have
\begin{align*}
0 &= \Phi_{c}(\gamma(s),\gamma(s))\\ &\leq
\Phi_{c}(\gamma(s),\gamma(t))+\Phi_{c}(\gamma(t),\gamma(s))\\
&\leq
A_{L+c}(\gamma|_{[s,t]})+\Phi_{c}(\gamma(t),\gamma(s))\\
&\leq 0
\end{align*}
in view of \eqref{aaa}. This implies that $\gamma$ is a static curve.

Suppose now that $\gamma:\R\to X$ is a static curve. Then $\gamma$ is a semistatic
curve with energy $c$ and given $s < t$ and $\epsilon>0$, there exists a curve
$\bar\gamma$ connecting $\gamma(t)$ to $\gamma(s)$ such that
$$ A_{L+c}\left(\gamma|_{[s,t]}\right)+ A_{L+c}(\bar\gamma)\leq
\epsilon. $$ Looking at the loop formed by $\gamma|_{[s,t]}$ and $\bar\gamma$, we
conclude that $h(\gamma(t),\gamma(t)) \leq 0$. But $h(x,x) \geq \Phi_c(x,x) = 0$ for
all $x\in X$, and hence $h(\gamma(t),\gamma(t)) = 0$. It follows from
\eqref{peierls} that $\gamma(t)\in \cA$, and thus $(\gamma(t),\dot{\gamma}(t))\in
\tilde\cA$, as we wanted to prove.
\end{proof}

\section{Minimizing optical hypersurfaces}
\label{minsur}

In this section, we show that many of the concepts that we presented in
Section~\ref{abm} do not really depend on the Lagrangian (or the Hamiltonian), but
can rather be formulated in the more general framework of optical hypersurfaces.

Let $\theta: T^*X\to X$ be the cotangent bundle of a closed manifold $X$, equipped
with the canonical symplectic form $\omega=d\lambda$, where $\lambda$ is the
Liouville 1-form. Let $\Sigma\subset T^*X$ be an optical hypersurface. Denote by
$\sigma$ its characteristic foliation. Recall that $\sigma$ is orientable and we
choose the orientation defined by the Hamiltonian vector field of any Hamiltonian
function which is fiberwise strictly convex and has $\Sigma$ as a regular level set.
Denote by $U_{\Sigma}$ the closed domain bounded by $\Sigma$.

\begin{defn}
An optical hypersurface $\Sigma$ is \emph{minimizing} if the interior
of $U_{\Sigma}$ does not contain a Lagrangian submanifold from
$\cL_0$, but any open neighborhood of $U_{\Sigma}$ does.
\end{defn}

\begin{rem}
\label{explain}
\begin{enumerate}

\item Theorem~\ref{corf1} ensures that we may replace ``Lagrangian
submanifold from $\cL_0$'' by ``Lagrangian section from $\cL_0$'', and
obtain precisely the same concept.

\item Suppose $\Sigma$ is a minimizing optical hypersurface, and pick any convex
superlinear Hamiltonian $H$ that has $\Sigma$ as a regular level set $H^{-1}(h)$.
Then, in view of \eqref{newway}, $h=c(H)$ is the Ma\~n\'e critical value of $H$.

\item If an optical hypersurface $\Sigma$ contains a Lagrangian submanifold
$\Lambda\in \cL_0$ (e.g., an exact Lagrangian section) then $\Sigma$ is minimizing.
Indeed, by Gromov's theorem \cite{gromov}, any two Lagrangian submanifolds
$\Lambda,K\in\cL_0$ must intersect, so the interior of $U_{\Sigma}$ cannot contain
an element of $\cL_0$.

The converse is certainly not true. However, Fathi explained \cite{F2} that if
$\Sigma$ is minimizing and every open neighborhood of $\Sigma$ contains a
$C^{\infty}$ exact Lagrangian section, then $\Sigma$ contains a $C^1$ exact
Lagrangian section. On the other hand, Fathi and Siconolfi \cite{FS} recently proved
that there always exists a $C^1$ function $f:X\to\R$ such that, first, $(x,df(x))\in
U_{\Sigma}$ for all $x\in X$ and, secondly, $(x,df(x))\in \Sigma$ if and only if
$x\in \cA$.
\end{enumerate}
\end{rem}

The last remark prompts the following question.

\begin{que}
Suppose an optical hypersurface $\Sigma$ contains an exact Lagrangian submanifold
$\Lambda\notin \cL_0$. Is $\Sigma$ minimizing?
\end{que}

\begin{rem}
The answer is ``Yes'' if $K$ admits a generating function quadratic at infinity (see
Definition~\ref{defgf}). Indeed, in this case, Theorem~\ref{lip} guarantees the
existence of a graph selector, and Theorem~\ref{corf1} shows that every neighborhood
of $U_{\Sigma}$ does contain a Lagrangian submanifold (even a section) from $\cL_0$.
\end{rem}

In the following, we are going to replace the concept of minimizing measure for a
convex Lagrangian $L$ by a notion that depends only on the foliation $\sigma$ of an
energy surface, and not on the particular choice of $L$. The appropriate notion is
that of \emph{foliation cycle} introduced by D. Sullivan in \cite{Su}. We briefly
review these ideas.

Let $M$ be a closed $n$--dimensional manifold and let $\Omega_{p}$ be the real
vector space of smooth $p$--forms on $M$. This vector space has a natural topology
which makes it a locally convex linear space. A continuous linear functional
$f:\Omega_{p}\to\R$ is called a $p$--current. Let ${\mathcal D}_p= (\Omega_p)^*$ be
the real vector space of all $p$--currents. With a natural topology, ${\mathcal
D}_p$ also becomes a locally convex linear space. Given a $p$--current $f$, we
define its boundary $\partial f$ as the $(p-1)$--current such that $\partial
f(\omega)=f(d\omega)$ for all $\omega\in \Omega_{p-1}$. Currents with zero boundary
are called cycles.

Among the set of all 1--currents, Sullivan considers a distinguished subset that he
calls \emph{foliation currents}. This subset is defined as follows. Given $x\in M$,
let $\delta_{x}: \Omega_{1}\to\R$ be the Dirac 1--current defined by
$\delta_{x}(\omega):= \omega_{x}(V(x))$. By definition, foliation currents are the
elements of the closed convex cone in $\mathcal{D}_1$ generated by all the Dirac
currents. A \emph{foliation cycle} is a foliation current whose boundary is zero.

Suppose now that $V$ is a non--vanishing vector field on $M$. Then $V$ defines a map
$\mu\mapsto f_{V,\mu}$ from measures to 1--currents, given by
$$ f_{V,\mu}(\omega) := \int_M \omega(V)\,d\mu. $$ Sullivan
\cite[Prop.~II.24]{Su} shows that this map yields continuous bijections between
\begin{enumerate}
\item nonnegative measures on $M$ and foliation currents; \item measures on $M$,
invariant under the flow of $V$, and foliation cycles.
\end{enumerate}

In our setting, $M$ is a minimizing optical hypersurface $\Sigma\subset T^*X$. Pick
some fiberwise convex, superlinear Hamiltonian $H$ such that $\Sigma= H^{-1}(h)$ is
a regular level set, and let $L$ be the corresponding Lagrangian. In view of
Remark~\ref{explain}, we have $h=c$. The following simple observation allows us to
translate the notion of minimizing measure into the languange of foliation cycles of
the characteristic foliation. Namely, if $(x,v)$ is a point in the critical energy
level $E^{-1}(c)\subset TX$ then
\begin{equation}
\label{coboundary}
L(x,v)+c=\lambda(d\ell(V(x,v))),
\end{equation}
where $\lambda$ is the Liouville form, $V$ the Euler--Lagrange vector field, and
$\ell$ the Legendre transform. Now, by \eqref{ergodic}, an invariant measure $\mu$
is minimizing if $\int_{TX} (L+c)\,d\mu =0$. We also know from \cite{DC} that
minimizing measures have their support contained in the energy level $E^{-1}(c)$.
Hence, the correct translation of the notion of minimizing measures into the
language of foliation cycles is the following.

\begin{defn}
Let $\Sigma$ be a minimizing optical hypersurface in $T^*X$, and $\sigma$ its
characteristic foliation. A foliation cycle $f$ of $\sigma$ is called
\emph{minimizing} if, and only if, $f(\lambda)=0$.
\end{defn}

In other words, minimizing foliation cycles are precisely those which can be
represented by measures $\ell_* \mu$ on $T^*X$, where $\mu$ is some minimizing
measure for some Hamiltonian $H$ with $\Sigma= H^{-1}(h)$. Observe also that, if we
have two Hamiltonians $H_1,H_2$ with the same regular level set $\Sigma$, and two
minimizing measures $\mu_1,\mu_2$ of $H_1,H_2$ representing the same foliation cycle
$f$, then the supports $\mu_1$ and $\mu_{2}$ will coincide. Hence it makes sense to
talk about the support of a foliation cycle $f$ of $\sigma$.

Now, the \emph{Mather set} of $\Sigma$ is defined as the closure of the union of the
supports of all minimizing foliation cycles; it coincides with the Mather set
$\tilde\cM^*$ in $T^*X$ of any convex superlinear Hamiltonian $H$ having $\Sigma$ as
regular level set.

In order to go further and define the Aubry set of $\Sigma$, we first have to
explain what a weak KAM solution should be in our setting. Given a point $(x,p)\in
\Sigma$, let $\Gamma^{\pm}{(x,p)}$ be the oriented positive (respectively, negative)
half of the leaf $\Gamma_{(x,p)}$ of $\sigma$ through $(x,p)$.

\begin{defn}
Let $\Sigma$ be a minimizing optical hypersurface in $T^*X$. A function
$u_{+}:X\to\R$ is called a \emph{positive weak KAM solution} of $\Sigma$ if the
following two conditions hold:
\begin{enumerate}
\item $u_+$ is Lipschitz, and $(x,du_+(x))\in U_{\Sigma}$ for almost every $x\in X$;
\item for every $x\in X$, there exists $(x,p)\in\Sigma$ such that, if $(y,p')$ is
any point in $\Gamma^{+}_{(x,p)}$, then
$$ u_{+}(y)-u_{+}(x)=\int_{\Gamma^{+}_{(x,p)}(y,p')}\lambda, $$ where
$\Gamma^{+}_{(x,p)}(y,p')$ is the oriented part of the leaf between $(x,p)$ and
$(y,p')$.
\end{enumerate}

Similarly, a function $u_-: X\to\R$ is called a \emph{negative weak KAM solution} of
$\Sigma$ if the following two conditions hold:
\begin{enumerate}
\item $u_-$ is Lipschitz, and $(x,du_-(x))\in U_{\Sigma}$ for almost every $x\in X$;
\item for every $x\in X$, there exists $(x,p)\in\Sigma$ such that, if $(y,p')$ is
any point in $\Gamma^-_{(x,p)}$, then
$$ u_-(x)-u_-(y)=\int_{\Gamma^-_{(x,p)}(y,p')}\lambda, $$ where
$\Gamma^-_{(x,p)}(y,p')$ is the oriented part of the leaf between $(y,p')$ and
$(x,p)$.
\end{enumerate}
\end{defn}

Again, \eqref{coboundary} shows that the sets $\cS_{\pm}= \cS_{\pm}(\Sigma)$ of
positive (respectively, negative) weak KAM solutions depend only on $\Sigma$ and not
on the particular choice of $H$ (or $L$). Setting
$$ \cI_{(u_-,u_+)} := \{ x\in X \mid u_-(x)=u_+(x) \} $$ for a pair of conjugate
functions, we see as before that the functions $u_{\pm}$ are differentiable on
$\cI_{(u_-,u_+)}$ with the same derivative. Therefore, the map $x\mapsto
du_-(x)=du_+(x)$ defines a set $\tilde\cI_{(u_-,u_+)}$ in $T^*X$ that contains the
Mather set of $\Sigma$. The \emph{Aubry set} of $\Sigma$ in $T^*X$ is then given by
$$ \tilde\cA^* = \tilde\cA^*(\Sigma) = \bigcap_{(u_{-},u_{+})}
\tilde{\cI}_{(u_{-},u_{+})}, $$ where the intersection is taken over all pairs
$(u_-,u_+)$ of conjugate functions.

Having defined the Aubry set, one would now like to study the dynamics on it and
single out a certain dynamically relevant set inside the Aubry set. For this, we
need the following general definition.

\begin{defn}
\label{scrdef}
Let $\phi_t$ be a continuous flow on a compact metric space $(X,d)$. Given
$\epsilon>0$ and $T>0$, a \emph{strong $(\epsilon,T)$--chain} joining $x$ and $y$ in
$X$ is a finite sequence $\{(x_{i},t_{i})_{i=1}^{n}\}\subset X\times\R$ such that
$x_1=x, x_{n+1}=y, t_i>T$ for all $i$, and $\sum_{i=1}^n d(\phi_{t_i}(x_i),x_{i+1})
< \epsilon$.

A point $x\in X$ is said to be \emph{strongly chain recurrent} if for all
$\epsilon>0$ and $T>0$, there exists a strong $(\epsilon,T)$--chain that begins and
ends in $x$.
\end{defn}

We denote by $\cR$ the set of all strong chain recurrent points. It contains the
nonwandering set\footnote{A point $x\in X$ is nonwandering if, and only if, for
every neighborhood $U$ of $x$ there exists a $T>1$ such that $\phi_T(U)\cap U\neq
\emptyset$; this implies that there are also arbitrarily large $T$ with that
property.}, but it is easy to give examples showing that it could be strictly
larger. The notion of strong chain recurrence strengthens the usual notion of chain
recurrence where one requires only $d(\phi_{t_i}(x_i),x_{i+1}) < \epsilon$ for every
single $i$. Strong chain recurrence was probably first considered by R. Easton in
\cite{E}.

Given an smooth orientable 1--dimensional foliation $\sigma$ on a closed manifold,
the strong chain recurrent set of $\sigma$ is the strong chain recurrent set of the
flow of any non--vanishing vector field $V$ tangent to $\sigma$. In the case where
$\sigma$ is the characteristic foliation of a hypersurface $\Sigma\subset T^*X$, we
denote by $\cR^*(\sigma)\subset \Sigma$ the strong chain recurrent set in $T^*X$,
and by $\cR(\sigma)\subset TX$ its preimage under the Legendre transform.

\begin{thm}
\label{new2}
Let $\Sigma$ be a minimizing hypersurface in $T^*X$, and let $K\subset \Sigma$ be an
exact Lagrangian submanifold (not necessarily in $\cL$). Then
$$ \cR^*(\sigma|_K)\subset \tilde\cA^*(\Sigma). $$ In particular, $\cR^*(\sigma|_K)$
is a Lipschitz graph over $X$.
\end{thm}

\begin{proof}
Choose a convex superlinear Hamiltonian $H$ which has $\Sigma$ as a regular level
set and let $L$ be its associated Lagrangian. For the proof, we will work on $TX$.
Endow $TX$ and $X$ with auxiliar Riemannian distances $d_{TX}$ and $d_{X}$ so that
the natural projection $\tau: TX \to X$ does not increase the distances. Consider
$(x,v)\in {\cR}$ and let $\gamma(t)=\tau(\phi_{t}(x,v))$. In view of
Proposition~\ref{aubdes} it suffices to show that the curve $\gamma$ is static. This
will imply that $\cR\subset \tilde\cA$.

Take $s\leq t$ and let $\xi:=\phi_{s}(x,v)$ and $\eta=\phi_{t}(x,v)$. We claim that,
for any $\epsilon>0$, there exists a strong $(\epsilon,1)$--chain that goes from
$\eta$ to $\xi$. To see this, let us start with a strong $(\delta,T)$--chain from
$x_1=\eta$ to $x_{n+1}=\eta$ where $T>1$ is large compared to $t-s$, and replace
$x_{n+1}$ by $\phi_{t_n-(t-s)}(x_n)$. If $\delta>0$ is chosen sufficiently small,
the point $\phi_{t_n-(t-s)}(x_n)$ lies in an $\epsilon$--neighborhood of $\xi$, and
we obtain a strong $(\epsilon,1)$--chain from $\eta$ to $\xi$. Let us call this
chain $\{(\eta_{i},t_{i})_{i=1}^{n}\}$ with $\eta_{1}=\eta$, $\eta_{n+1}=\xi$,
$t_{i}>1$ and $\sum_{i=1}^{n}d_{TX}(\phi_{t_{i}}(\eta_{i}),\eta_{i+1})<\epsilon$.
Set $p_{i}:=\tau(\eta_{i})$ and $q_{i}:=\tau(\phi_{t_{i}}(\eta_{i}))$.

Recall that $\Phi_{c}$ satisfies
$$ \Phi_c(x,y)\leq \Phi_c(x,z)+\Phi_c(z,y) $$ and
$\Phi_{c}(x,x)=0$ for all $x,y,z\in X$. Hence
$$ \Phi_{c}(p_{1},p_{n+1})\leq \Phi_{c}(p_1,q_1)+
\Phi_{c}(q_1,p_2)+\cdots+ \Phi_{c}(p_n,q_n)+ \Phi_{c}(q_n,p_{n+1}). $$ Given $p$ and
$q$ in $X$, let $\gamma:[0,d_X (p,q)]\to X$ be a unit speed minimizing geodesic
connecting $p$ to $q$. We have
$$ \Phi_{c}(p,q)\leq \int_{0}^{d_X
(p,q)}(L+c)(\gamma(t),\dot{\gamma}(t))\,dt\leq \kappa_1 d_X (p,q) $$ where $\kappa_1
:=\max\{|(L+c)(x,v)| \mid (x,v)\in TX \text{ and } |v|=1\}$. Thus
\begin{equation}
\label{esti1}
\sum_i \Phi_c(q_i,p_{i+1}) \leq \kappa_1 \sum_i d_X(q_{i},p_{i+1}) \leq \kappa_1
\epsilon.
\end{equation}
Using \eqref{coboundary} and the fact that $K$ is exact, we have
\begin{equation}
\label{esti2}
\Phi_c(p_i,q_i) \leq A_{L+c}(\tau\circ\phi_{t_i}(\eta_i)|_{[0,t_i]}) =
g(\phi_{t_i}(\eta_i)) - g(\eta_i),
\end{equation}
where $g:\ell^{-1}(K)\to\R$ is a smooth function such that
$d(g\circ\ell^{-1}|_K)=\lambda|_K$. Combining \eqref{esti1} and\eqref{esti2}, we
obtain
$$ \Phi_{c}(p_{1},p_{n+1}) \leq \sum_i \Phi_c(p_i,q_i) + \Phi_c(q_i,p_{i+1}) \leq
\kappa_1 \epsilon + \kappa_2 \epsilon + g(\xi)-g(\eta), $$ where $\kappa_2$ is a
Lipschitz constant for $g$. On the other hand,
$$ \Phi_{c} (p_{n+1},p_1) \leq A_{L+c}(\gamma|_{[s,t]}) =
g(\eta)-g(\xi). $$ Therefore
$$ 0 = \Phi_c(p_1,p_1) \leq \Phi_c(p_1,p_{n+1}) + \Phi_c(p_{n+1},p_1) \leq
(\kappa_1 + \kappa_2)\epsilon.  $$ Since $\epsilon$ is arbitrary, we have
$$\Phi_c(p_1,p_{n+1}) + \Phi_c(p_{n+1},p_1)= 0. $$ Using the triangle
inequality for $\Phi_c$ as in the proof of Proposition~\ref{aubdes}, we see that
$\gamma$ is a static curve.
\end{proof}

\section{Non--removable intersections in the convex case}
\label{nonrem0}

Let $\Sigma\subset T^*X$ be an optical hypersurface bounding a domain $U_{\Sigma}$,
and $\tilde\cA^*=\tilde\cA^*(\Sigma)$ its Aubry set in $T^*X$. The following theorem
is the main result that combines non--removable intersections with Aubry--Mather
theory. Its statement was pointed out to us by A. Fathi, who also explained how
results from Nonsmooth Analysis could be used to simplify its proof.

\begin{thm}
\label{new1}
Let $\Sigma$ be a minimizing optical hypersurface such that $U_{\Sigma}$ contains a
Lagrangian submanifold $\Lambda\in \cL_0$. Then $$ \tilde\cA^* \subset \Lambda\cap
\Sigma. $$
\end{thm}

\begin{proof}
Let $u:X\to \R$ be a graph selector associated to $\Lambda$ (see Theorem~\ref{lip}).
The function $u$ is Lipschitz and $(x,du(x))\in \Lambda$ for almost every $x\in X$.
By Remark~\ref{smaller} and Theorem~\ref{conjugate}, there exists a pair of
conjugate functions $(u_{-},u_{+})$ with $u_{+}\leq u\leq u_{-}$. At any point $x\in
\cI_{(u_{-},u_{+})}$, the three functions are differentiable with the same
derivative. Hence $du(x)$ exists for each $x\in \cI_{(u_{-},u_{+})}$ and satisfies
$(x,du(x))\in\Sigma$.

We will now show that for every $x\in \cI_{(u_{-},u_{+})}$ we have $(x,du(x))\in
\Lambda$. This is the main difficulty since, a priori, we only know that this is
true for almost every $x$. For each $x\in X$, let $C_{x}(\Lambda)$ denote the convex
hull of $\Lambda\cap T^*_{x}X$. Since $\Lambda\cap T^*_{x}X$ is compact,
$C_{x}(\Lambda)$ is also compact by Carath\'eodory's theorem. Let
$C(\Lambda)=\cup_{x\in X}C_{x}(\Lambda)$. Since $C(\Lambda)$ is compact, a result
from Nonsmooth Analysis (cf. \cite[pp.~62--63]{Clarke} and \cite[Prop.~8.4]{FM})
ensures that for any point $x$ of differentiability of $u$ we have $(x,du(x))\in
C(\Lambda)$. But since $(x,du(x))\in \Sigma$, $\Sigma\cap T_{x}^{*}X$ is strictly
convex and $\Lambda\subset U_{\Sigma}$, the point $(x,du(x))$ is an extreme point of
$C_{x}(\Lambda)$. But any extreme point in the convex hull must belong to
$\Lambda\cap T^*_xX$ and thus $(x,du(x))\in \Lambda$.

Since $\tilde\cA$ is contained in $\tilde\cI_{(u_{-},u_{+})}$, the theorem follows.
\end{proof}

This result can be applied to boundary rigidity. The following result is a
generalization of Theorem~\ref{brigid}.

\begin{thm}
\label{new3}
Let $\Lambda\in \cL$ be a Lagrangian submanifold lying in an optical hypersurface
$\Sigma$. Assume that $\sigma|_{\Lambda}$ is strongly chain recurrent. Let $K\in
\cL$ be any Lagrangian submanifold lying in $U_{\Sigma}$ with the same Liouville
class as $\Lambda$. Then $K=\Lambda$.
\end{thm}

\begin{proof}
Since the multidimensional Birkhoff theorem is valid if $\sigma|_{\Lambda}$ is chain
recurrent \cite[Prop.~1.2]{birkhoff}, we may, as in Section~\ref{rproof}, apply a
symplectic shift and assume that $\Lambda= \cO\subset T^*X$. The shifted
hypersurface obtained from $\Sigma$ is minimizing since it contains $\cO$ (see
Remark~\ref{explain}). We then know from Theorem \ref{new2} that $\cO \subset
\tilde\cA^*$. Since the natural projection $\theta|_{\tilde\cA^*}: \tilde\cA^* \to
\cA$ is a homeomorphism \cite[Prop.~5.2.8]{F}, we have $\tilde\cA^* = \cO$.

Now pick a graph selector $\Phi$ of $K$; see Theorem~\ref{lip}. As in the proof of
Theorem~\ref{new1}, it will be differentiable at every point in
$\theta(\tilde\cA^*)=X$ with zero derivative, i.e., $K$ coincides with the zero
section $\cO= \Lambda$.
\end{proof}

\begin{figure}
\begin{center}
\includegraphics[width=\textwidth]{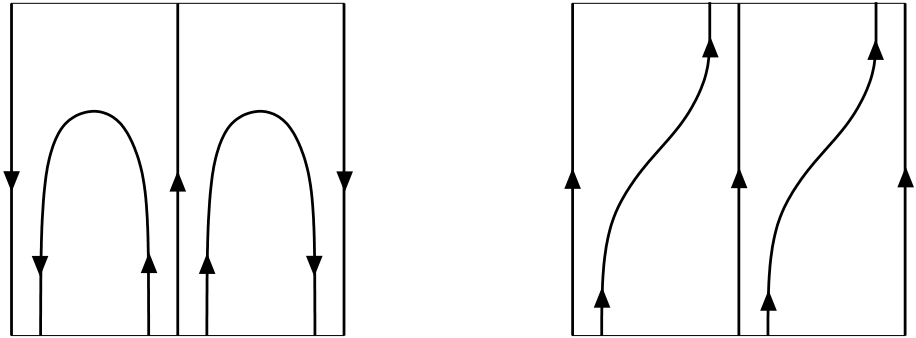}
\end{center}
\caption{The dynamics on $\cO$ in Theorem~\ref{twocycles} (left) and
Example~\ref{twotwo} (right)}
\label{pictt}
\end{figure}

\begin{proof}[Proof of Theorem~\ref{twocycles}]
Recall that we deal with the zero section $\cO$ of $T^*\T^2$ lying inside the
optical hypersurface
$$\Sigma = \{(y_1-\sin x_1)^2 + (y_2-\cos x_1)^2= 1\}.$$ The restriction
$\sigma|_{\cO}$ of the characteristic foliation is a Reeb foliation; see
Figure~\ref{pictt}. Denote by $Z$ the union of the two limit cycles. Note that $Z$
is the strong chain recurrent set of $\sigma|_{\mathcal O}$, and so, by
Theorem~\ref{new2}, we have $Z \subset \tilde\cA^*$. Since $\Sigma$ contains the
zero section $\cO$ it is minimizing in view of Remark~\ref{explain}. Applying
Theorem~\ref{new1}, we see that $Z \subset \Lambda \cap \Sigma$ for any Lagrangian
submanifold $\Lambda \in \cL _0$ lying in $U_{\Sigma}$. This completes the proof.
\end{proof}

As a by--product of the proof, we get the following explicit description of both the
Aubry and the Mather sets of $\Sigma$ in this situation:
$$\tilde\cA^*= \tilde\cM^* = Z.$$
Indeed, we have seen in Example~\ref{notrigid} that the graph of $df$ with
$f(x_{1},x_{2})=-\cos x_{1}$ intersects $\Sigma$ exactly along $Z$. Hence, by
Theorem~\ref{new1}, we obtain $\tilde\cA^*  \subset Z$. Together with the opposite
inclusion established in the proof above, this yields $Z = \tilde\cA^* $.
Furthermore, each of the limit cycles in $Z$ is a foliation cycle. It vanishes on
the Liouville form since $\lambda|_{\cO} = 0$. Hence $\tilde\cM^* = Z$.

\begin{ex}
\label{twotwo}
Let $f:S^1\to S^1$ be a diffeomorphism with only two fixed points such that the
fixed points are neither attractors nor repellers.  Let $X$ be the unit norm vector
field on $\T^2$ obtained by suspending $f$. Write
$$ X(x_1,x_2)=(a_{1}(x_1,x_2),a_{2}(x_1,x_2))$$ and let $H$ be the
Hamiltonian
$$
H(x_1,x_2,y_1,y_2)=(y_{1}-a_{1}(x_1,x_2))^{2}+(y_{2}-a_{2}(x_1,x_2))^{2}. $$
Consider $\Sigma = \{ H= 1 \} \subset T^* \T^2$. Since $\Sigma$ contains the zero
section $\cO$ it is minimizing in view of Remark~\ref{explain}. If we identify $\cO$
with $\T^2$ then $X$ is tangent to the characteristic foliation $\sigma|_{\cO}$.
Note that $\sigma|_{\cO} $ is strongly chain recurrent, hence $\cO$ is boundary
rigid by Theorem~\ref{new3}.

In this example, one can also describe the Aubry and the Mather sets of $\Sigma$
explicitly. Indeed, Theorems \ref{new2} and \ref{new1} yield  $\tilde\cA^* \subset
\cO$ and $\cO \subset \tilde\cA^*$, respectively, so $\tilde\cA^* = \cO$. The Mather
set is strictly contained in the Aubry set. Indeed, the support of any invariant
measure of the vector field $X$ lies in the union $Z$ of the limit cycles, each of
which vanishes on the Liouville form; hence $\tilde\cM^* = Z$.
\end{ex}

\section{ Constructing Lagrangian sections---proof of Theorem~\ref{corf1}}
\label{2.3}

Suppose $U\subset T^* X$ be an open, fiberwise convex set. We want to prove that
every class $a\in \sh{U}$ can be represented by a Lagrangian \emph{section} of the
cotangent bundle. Indeed, this an immediate consequence of the following theorem (a
more general version of which was proved independently in \cite[Appendix]{FM}). Let
us denote the fiberwise convex hull of a set $S\subset T^* X$ by $\conv{S}$.

\begin{thm}
\label{fund}
Given a Lagrangian submanifold $\Lambda\in \mathcal{L}$, the fiberwise convex hull
$\conv{W}$ of any neighborhood $W$ of $\Lambda$ contains a Lagrangian section
$\Lambda_0 \in \mathcal{L}$ with $a_{\Lambda_0} = a_{\Lambda}$.
\end{thm}

\begin{proof}
We may assume that $\Lambda$ is an exact Lagrangian submanifold, simply by applying
the symplectic shift $(x,y)\mapsto (x,y-\alpha(x))$ where $\alpha$ is the closed
1--form on $X$ representing the Liouville class $a_{\Lambda}$.

Let $\Phi: X\to\R$ be a graph selector associated to $\Lambda$ as described in
Theorem~\ref{lip}; namely, $\Phi$ is Lipschitz continuous, smooth on an open subset
$X_0 \subset X$ of full measure, and satisfies
\begin{equation}
\label{grsel}
\text{graph}(d\Phi|_{X_0}) \subset \Lambda.
\end{equation}
The proof of Theorem~\ref{fund} is divided into two steps.

\medskip
\noindent
\emph{Smoothing:} We are going to regularize the Lipschitz function
$\Phi$ by a convolution argument, similar to the proof of Prop.~7 in
\cite{CIPP1}. For this, we embed $X$ into some Euclidean space
$\R^N$. Denote by $V_r$ the $r$--neighborhood of $X$ in $\R^N$ where
$r>0$ is chosen small enough so that the orthogonal projection $\pi
:V_r \to X$ is well defined. We extend $\Phi: X\to\R$ to a function
$\bar\Phi: V_r\to\R$ by setting
$$ \bar\Phi := \Phi \circ \pi . $$

For each $s\in (0,r/2)$, we pick a smooth cut--off function $u:
[0,\infty)\to [0,\infty)$ with support in $[0,s]$ such that $u$ is
constant near 0 and satisfies
$$ \int_{\R^N} u(|z|)dz = 1 . $$ Define the function $\bar\Psi: V_s
\to \R$ as the convolution
$$\bar\Psi(z) := (\bar\Phi \ast u)(z) = \int_{R^N} \bar\Phi(y)
u(|z-y|) dy .$$

Since $\bar\Phi$ is Lipschitz continuous, it is differentiable almost
everywhere and weakly differentiable. Therefore, $\bar\Psi$ is a
smooth function on $V_s$ with
\begin{align*}
d\bar\Psi(z) &= \int_{\R^N} \bar\Phi(y) d_z u(|z-y|) dy \\ &= -
 \int_{\R^N} \bar\Phi(y) d_y u(|z-y|) dy \\ &= \int_{\R^N}
 d\bar\Phi(y) u(|z-y|) dy.
\end{align*}
Denote by
$$ \Psi := \bar\Psi|_X $$ the restriction of $\bar\Psi$ to $X$, and
let $B_s(x) \subset V_s \subset \R^N$ be the open ball of radius $s$
centered at $x\in X$. Because $X_0$ has full measure in $X$, we
conclude that
\begin{equation}
\label{psi} d\Psi(x) = \int_{\pi^{-1}(X_0) \cap B_s(x)}
d\bar\Phi(y)|_{T_xX} u(|x-y|)dy.
\end{equation}
Note that, for this formula to make sense, we identify each
$T_y\R^N$ (where $y\in \R^N$) with $\R^N$, and each $T_xX$
(where $x\in X$) with a linear subspace of $\R^N$.

\medskip
\noindent \emph{Analysing formula \eqref{psi}:} For each $x
\in X$, we write
$$ P_x: T_x\R^N \cong \R^N \to T_xX $$ for the orthogonal projection. Write
$|\cdot|$ for the Euclidean norm on $\R^N$ and
$|\cdot|^*$ for the dual norm on $(\R^N)^*$.
Introduce a distance function on $T^*X$
by setting
\begin{equation}
\label{dist} \text{dist}((x,\xi),(y,\eta)) := |x-y| +
|\xi\circ P_x - \eta\circ P_y |^*.
\end{equation}

For $x\in X$, we define the set
$$ \mathcal{G}_s(x) := \{(x, d\bar\Phi(y)|_{T_xX})) \mid y\in \pi^{-1}(X_0) \cap
B_s(x) \} \subset T^*X. $$ For a subset $Z \subset T^*X$, we
denote by $W_{\epsilon}(Z)$ the $\epsilon$--neighborhood of
$Z$ with respect to the distance defined in \eqref{dist}.

\begin{lemma}
\label{anal}
For every $\epsilon > 0$ there is an $s>0$ such that
$$ \mathcal{G}_s(x) \subset W_{\epsilon/2}(\emph{graph}(d\Phi|_{X_0})) $$
for each $x \in X$.
\end{lemma}

\begin{proof}
Pick any point
$$ \eta_1 = (x,d\bar\Phi(y)|_{T_xX}) \in \mathcal{G}_s(x)$$
with $x\in X$ and $y\in \pi^{-1}(X_0)\cap B_s(x)$. We will
show that the distance between $\eta_1$ and
$$ \eta_2 := (\pi(y),d\Phi(\pi(y))) \in \text{graph}(d\Phi|_{X_0})
$$ becomes as small as we wish when $s\to 0$ uniformly in $x$ and $y$.
Indeed, denote by $c>0$ the Lipschitz constant of $\Phi$ with respect to the induced
distance on $X \subset \R^{N}$. Let $Q_y$ be the differential of the projection
$\pi$ at $y$; we consider $Q_y$ as an endomorphism of $\R^N$. Finally, write $\Vert
\cdot \Vert$ for the operator norm on $\text{End}(\R^N)$. Now we can estimate
\begin{align*}
\text{dist}(\eta_1,\eta_2) &= |x-\pi(y)| +
|d\bar\Phi(y)|_{T_xX}\circ P_x - d\Phi(\pi(y))\circ P_{\pi(y)}|^*
\\ &= |x-\pi(y)| + |d\Phi(\pi(y))\circ
Q_y \circ P_x - d\Phi(\pi(y))\circ P_{\pi(y)}|^*
\\ &\leq |x-y| +|y-\pi(y)| + c \Vert
Q_y \circ P_x - P_{\pi(y)} \Vert.
\end{align*}
Note that $|x-y|+|y-\pi(y)| \leq 2s \to 0$
as $s\to 0$. It remains to handle the term
$ \Vert
Q_y \circ P_x - P_{\pi(y)} \Vert$.
Using that
$\Vert P_x\Vert = \Vert P_{\pi(y)}\Vert =1$
we obtain
\begin{align*}
\Vert Q_y \circ P_x - P_{\pi(y)}\Vert &=
\Vert Q_y \circ P_x - P_{\pi(y)}\circ P_x + P_{\pi(y)}\circ P_x
-P_{\pi(y)} \circ P_{\pi(y)} \Vert \\ &\leq
\Vert Q_y - P_{\pi(y)}\Vert + \Vert P_x - P_{\pi(y)}\Vert \to 0
\end{align*}
as $s \to 0$, and the convergence is uniform in
$x\in X$ and $y\in B_s(x)$.
This finishes the proof of
Lemma~\ref{anal}.
\end{proof}

Now we can readily prove Theorem~\ref{fund}. Namely, given any $\epsilon > 0$, we
choose $s$ as in Lemma~\ref{anal}. Then \eqref{psi}, Lemma~\ref{anal}, and
\eqref{grsel} imply that
$$ (x,d\Psi(x)) \in \conv{W_{\epsilon/2}(\mathcal{G}_s(x))} \subset
\conv{W_{\epsilon}(\text{graph}(d\Phi|_{X_0}))} \subset \conv{W_{\epsilon}(\Lambda)}
$$ for each $x \in X$. Thus the Lagrangian section $\Lambda_0 := \text{graph}(d\Psi)$
satisfies
$$ \Lambda_0 \subset \conv{W_{\epsilon}(\Lambda)}. $$ Since $\epsilon >0$ was
arbitrary the proof of Theorem~\ref{fund} is completed.
\end{proof}

\section{Boundary rigidity in general symplectic manifolds}

The boundary rigidity phenomenon can be naturally formulated in the following more
general context. Let $(M,\omega)$ be a compact symplectic manifold with non--empty
boundary, and let $\Lambda \subset {\partial M}$ be a closed Lagrangian submanifold.
Denote by $\cL_0$ the space of all Lagrangian submanifolds of $M$ which are exact
Lagrangian isotopic to $\Lambda$ (see for instance \cite{Pbook} for the definition
of exact Lagrangian isotopies in symplectic manifolds). We say that $\Lambda$ is
\emph{boundary rigid} if $\cL_0 = \{\Lambda\}$, and \emph{weakly boundary rigid} if
every $K \in \cL_0$ is contained in $\partial M$.

Theorem~\ref{brigid} already provides a class of examples of
boundary rigid Lagrangian submanifolds. It would be
interesting to investigate boundary rigidity in other
symplectic manifolds as well.

\begin{ex}[A toy example]
Let $M=D^2$ be the 2--disc endowed with some area form. Then the circle $\Lambda =
\partial{D^2}$ is boundary rigid. Indeed, every circle $K \in \cL_0$ must enclose
the same area as $\Lambda$, and hence $K=\Lambda$.
\end{ex}

Can one generalize this example to higher dimensions? For
instance, let $M$ be the Euclidean ball
$$ \{p_1^2 + q_1^2 +p_2^2+q_2^2 \leq 2\} $$
in the standard symplectic vector space $\R^4$. Consider the
split torus
\begin{equation}
\label{split} \Lambda := \{p_1^2 +q_1^2 = 1,\;p_2^2+q_2^2 = 1\}
\subset \partial{M}.
\end{equation}
One can show that $\Lambda$ admits a nontrivial exact Lagrangian isotopy inside
$\partial M$ and hence is \emph{not} boundary rigid.

\begin{que}
\label{qu0}
Is $\Lambda$ weakly boundary rigid? What happens with Lagrangian tori contained in
general ellipsoids in $\R^{2n}$?
\end{que}

Further, it would be interesting to extend the study of non--removable
intersections, both between Lagrangian submanifolds, and between a Lagrangian
submanifold and a hypersurface (see Sections~\ref{nrint} and \ref{nonrem0}), to more
general symplectic manifolds.

An interesting playground for this problem is given by tori in ellipsoids as in
Question~\ref{qu0}. For instance, let $\Lambda$ be the split torus given by
\eqref{split}, and $K$ any exact Lagrangian deformation of $\Lambda$ which lies in
the closed ball $M$ of radius $\sqrt{2}$. Just recently, Y. Eliashberg outlined a
beautiful argument based on symplectic field theory which suggests that $K$ must
intersect the boundary of the ball; later, F. Schlenk proposed a simpler approach
using symplectic capacities, based on \cite{vit}. Applying Sullivan's theorem as in
the proof of Theorem~\ref{nonrem1}, we conclude that $K \cap \partial M$ must
contain a closed orbit of the characteristic foliation of $\partial M$.

We refer to \cite{Bir} for further discussion on symplectic
intersections.

\bigskip \noindent \emph{Acknowledgement}. Our interest in the boundary rigidity
phenomenon and, in particular, non--removable intersections was motivated by
numerous fruitful discussions with Paul Biran. Discussions with Albert Fathi were
very helpful; in particular, he explained to us how to use results from Nonsmooth
Analysis in the proof of Theorem \ref{new1}. We are also grateful to Misha Bialy and
Yasha Eliashberg for useful conversations.

Parts of this paper were written during joint visits to Tel Aviv, Cambridge, and
Z\"urich. These visits were supported by the Hermann Minkowski Center for Geometry
at Tel Aviv University, the University of Cambridge, and the Forschungsinstitut
f\"ur Mathematik at ETH Z\"urich. We thank these institutions for their support, as
well as Marc Burger, Kostya Khanin, Dietmar Salamon, and Edi Zehnder for their warm
hospitality.

\end{document}